\documentclass[12pt]{article}

\usepackage{amsfonts}
\usepackage{latexsym}
\usepackage{amssymb}
\usepackage{amsmath}

\def\p{{\bf P}} 
\def\E{{\bf E}} 
\def\u{{\bf u}} 
\def\v{{\bf v}}  
\def\unit{{\bf 1}} 
\def\f{{\bf f}} 
\def\s{{\bf s}} 
\def\x{{\bf x}}
\def\y{{\bf y}}
\def\pp{{\bf p}} 
\def\e{{\rm e}}
\def\d{{\rm d}} 
\def\ed{{\stackrel{{\cal D}}{=}}} 
\def\R{{\mathbb R}} 
\def\N{{\mathbb N}} 
\def\h{{\mathbf h}} 
\def\z{{\mathbf z}} 
\def\supp{\textrm{supp}} 
\def\k{{\mathbf k}} 
\def\K{{\mathbf K}} 
\def\X{{\mathbf X}} 
\def\Rd{\R^d}

\newtheorem{lemma}{Lemma}
\newtheorem{theorem}{Theorem}

\newenvironment{proof*}{\smallskip \noindent {\bf Proof.} \ }
{\nopagebreak \hspace*{1pt} \vspace{-2pt} \hfill\rule{6pt}{6pt} \medskip}

\newenvironment{theorem*}{\medskip \noindent{\bf Theorem.} \it }{\medskip} 

\begin{document}

\begin{center}
{\Large \textbf{Critical Multitype Branching Systems: Extinction Results}} \\

\vspace{1cm}

\textsc{P\'eter Kevei}
\footnote{Supported by the Analysis and Stochastics Research Group of the Hungarian
Academy of Sciences.}   \\

Centro de Investigaci\'on en Matem\'aticas \\ 
Jalisco S/N, Valenciana, Guanajuato, GTO 36240, Mexico\\
e-mail: \texttt{kevei@cimat.mx} \\

\bigskip

\textsc{Jos\'e Alfredo L\'opez Mimbela}

Centro de Investigaci\'on en Matem\'aticas \\
Jalisco S/N, Valenciana, Guanajuato, GTO 36240, Mexico\\
e-mail: \texttt{jalfredo@cimat.mx}
\end{center}

\bigskip

\begin{abstract}
We consider a critical branching particle system in  $\R^d$, composed of individuals
of a finite number of types $i\in\{1,\ldots,K\}$. Each individual of type $i$ moves independently according
to a symmetric $\alpha_i$-stable motion. We assume that the particle lifetimes and offspring
distributions are type-dependent. Under the usual independence assumptions in branching systems, we prove extinction theorems in the following cases:
(1) all the particle lifetimes have finite mean, or
(2) there is  a type whose lifetime distribution has
heavy tail, and the other lifetimes have finite mean.
We get a more complex dynamics by assuming in case (2) that the most mobile particle
type corresponds to a finite-mean lifetime:
in this case, local extinction of the population is determined by an interaction of
the parameters (offspring variability, mobility, longevity) of the long-living
type and those of the most mobile type.
The proofs are based on  a precise analysis of the occupation times of a related
Markov renewal process, which is of independent interest.

\textit{Keywords:} Critical branching particle system; Extinction; Markov renewal process.

\textit{AMS Subject Classification:} MSC 60J80, MSC 60K15.
\end{abstract}

\section{Introduction}

In critical branching and migrating populations, mobility of individuals counteracts the tendency to asymptotic local extinction caused by the clumping effect of the branching. In fact, convergence to a non-trivial equilibrium may occur in a spatially distributed population whose members perform migration and reproduction, even if the branching is critical, provided that the mobility of individuals is strong enough. This behavior has been investigated in
several branching models, including branching random walks   \cite{Kallenberg,Matthes}, Markov branching systems (both with monotype \cite{GW}  and multitype \cite{FV,GRW,LW} branching), and age-dependent  branching systems \cite{VW}.

In \cite{VW} Vatutin and Wakolbinger
investigated a monotype branching model in Euclidean space $\Rd$, in which each
particle moves according to a symmetric $\alpha$-stable motion, and at the end of its
 lifetime it leaves   at its death site  a random number of offsprings, with critical offspring generating function
$
f(s) = s + \frac{1}{2} ( 1 - s)^{1 + \beta}
$,
$\beta \in (0,1]$.
It turned out that, if the initial population is Poisson with uniform intensity and
the particle lifetime distribution has finite mean, such process suffers local extinction if $d \leq \alpha/\beta$, while
for $ d > \alpha/\beta$ the system is persistent, i.e.~preserves its intensity in the large
time limit. This result is consistent with the intuitive meaning of the population parameters: the exponent
$\alpha>0$ is the mobility parameter of individuals in the sense that a smaller $\alpha$ means a more mobile migration (i.e. more spreading out of  particles) which is clearly in favor of persistence;
$\beta$ is the offspring variability parameter, meaning that a smaller $\beta$
causes a stronger
clustering effect in the population, which
favors local extinction due to criticality of the branching. If the lifetime distribution has a power tail
$t^{-\gamma}$ for some $\gamma \in (0,1]$, then the critical dimension is
$\alpha \gamma/ \beta$. Again, it is intuitively clear that long lifetimes (i.e.~small
$\gamma$) enhance the  spreading out of individuals. However, Vatutin and Wakolbinger discovered that, in contrast with the case of finite-mean lifetimes, if the lifetimes have a general distribution of the above sort, the ``critical'' dimension does not
necessarily pertains to the local extinction regime: when $d=\alpha\gamma/\beta$ persistence of the population is not excluded.

Our aim in  the present paper is to get a better understanding  about  how population characteristics such as mobility, offspring variability, and longevity of individuals determine the asymptotic local extinction of branching  populations. In order to attain this we deal with a multitype system, where the most mobile
migration (corresponding to the smallest $\alpha$) and the life-time distribution with the heaviest tail, may correspond to different particle types.  More precisely, we consider a branching population living in  $\R^d$, constituted of particles of different types $i\in\K:=\{1,\ldots,K\}$.
 Each particle of type $i$ moves according to a symmetric
$\alpha_i$-stable motion
until the end of its random lifetime, which has a non-arithmetic distribution function ${\Gamma}_i$.
Then it branches according to a multitype offspring distribution with generating function
$f_i(\s)$, $\s\in[0,1]^K$, $i\in\K$. The descendants appear where the parent individual died,
and evolve independently in the same manner. The movements, lifetimes and branchings of particles are assumed to be independent; the only dependency in the system is that the offsprings start where the parent particle
died. In addition, we assume that the process starts off at time 0 from
a Poisson random population, with a prescribed intensity measure, and that all
particles at time 0 have  age 0.
Let $M = (m_{i,j})_{i,j=1}^K$
denote the mean matrix of the multitype branching law, that is
$$
m_{i,j}= \frac{\partial f_i}{\partial x_j}(\unit),
$$
where $\unit = (1,1, \ldots, 1) \in \R^K$.
We assume that $\f(\s) = (f_1(\s), \ldots, f_K(\s)) \ne M \s$,
and that $M$ is an ergodic stochastic matrix.
This implies that the branching is critical, i.e. the largest
eigenvalue of $M$ is 1.

For the system described above, here we investigate  parameter configurations under which the population becomes locally extinct in the large time run.
We deal first with the case when all particle lifetimes have finite mean and prove that the  process suffers local extinction if $d < \alpha/\beta$, where
the mobility parameter $\alpha= \min_{1 \leq i \leq K} \alpha_i$ is the same as
in the Markovian case \cite{LW}, and the offspring variability parameter
$\beta \in (0,1]$ is determined by
$$
x - \langle \v, \unit - \f(\unit - \u \, x) \rangle \sim x^{1+\beta} L(x)\mbox{ as } x\to0,
$$
where
 $\v$  denotes the (normalized) left
eigenvector of the matrix $M$ corresponding to the eigenvalue 1, and
 $L$ is slowly varying at 0 in the sense that
$\lim_{x \to 0}L(\lambda x)/ L(x) = 1$
for every $\lambda > 0$.
In a way, this case is similar to the one with  exponentially distributed lifetimes.

Next we assume that exactly one particle type is long-living, i.e. its lifetime distribution has a
power tail decay $t^{-\gamma}$, $\gamma \in (0,1]$, while the other lifetime types have distributions with tails decaying not slower than $A \, t^{-\eta}$ for some $\eta>1$, $A > 0$.
We consider two scenarios. In the first one we assume that the most mobile particle type is, at the same time,
long-living, and we prove that extinction holds when $d < \alpha \gamma/ \beta$.
Then we proceed with the
 most interesting case:
the most mobile particle type corresponds to a finite-mean lifetime.
In this scenario, it turns out that local extinction of the population is determined by a complex interaction of
the parameters (offspring variability, mobility, longevity) of the long-living type and those of the most mobile type. 
Assuming without loss of generality that type 1 is the long-living type, we prove that the systems suffers local extinction provided that $d<d_+$, where
$$
d_+ = \frac{ \gamma }
{\frac{( \beta + 1) \gamma}{\alpha} - \frac{1}{\alpha_1}}.
$$
The positive number $\gamma\alpha_1$ can be considered as the ``effective mobility'' parameter of the long-living type. If $\gamma\alpha_1$ is very close to $\alpha$ (so that $\gamma/\alpha$ and $1/\alpha_1$  are approximately the same), then
$d_+$ is also close to $\alpha_1 \gamma / \beta$ and to $\alpha/\beta$.
Moreover,
for fixed $\alpha, \alpha_1$ and $\gamma$, the parameter $d_+$ considered as a function of $\beta$,
is decreasing, which is consistent with previous known results.

The proofs of our results rely on a precise asymptotic analysis for the occupation times in the branching particle system. Some of our techniques  combine  parts of the approaches used in \cite{FV} and \cite{VW} adapted to our model, however the adaptation to our case  is far from being straightforward. In Section \ref{famtreeanal} we  provide a family
tree analysis which allows us to compare the occupation times of the particle
system with the occupation times of an auxiliary Markov renewal process.

Then, in Section \ref{markovren}, we carry out the asymptotic analysis mentioned above by investigating the occupation times of all types in the auxiliary renewal process, as well as the asymptotic number of renewals in large time-intervals.
This is the mathematical core of the paper, and we think it is interesting on its own right.
To achieve this, we need to control the tail decay of the renewal times of all types simultaneously, which we were able to do assuming that
there is only one long-living particle type. Therefore, in its present form our approach is not yet applicable to treat a general model with arbitrary lifetime distributions.

Finally, in Section \ref{extinction} we give the extinction results
in our various different setups.
Let us remark that, when the particle lifetimes have finite mean and the spatial dimension is small, local extinction of the population can be proved without the occupation times analysis; in this case a simple estimation yields the result, see the proof of Theorem \ref{finmean}.
In contrast, the occupation time analysis is needed to treat the case of long-living particle types.

\section{Family tree analysis} \label{famtreeanal}

 Following \cite{FV} (p.~553--558) we introduce the following auxiliary process.
Consider a Markov renewal process with values in $\K$,
where in type $i$ the process
spends time according to a non-lattice distribution $\Gamma_i$ (whose distribution function we denote again by $\Gamma_i$), such that $\Gamma_i(0)=0$, and then
jumps to type $j$ with probability
$m_{i,j}.$ We write $\mu_i=\int_0^\infty x \,\Gamma_i (\d x)$ for
the mean of the $i^{\textrm{th}}$ lifetime, which can be infinite.
Let $\overline t_j(t)$ be the time that the process spends at state $j$ up to time
$t$. Put $\overline r_{i,j} (t,a) = \p_i \left\{ \bar t_j(t) \geq a \right\}$,
$i,j \in \K$, where $\p_i$ stands for the probability when the process
starts in type $i$. We aim at finding an upper bound for the probabilities
$\overline r_{i,j}$.

First we show the connection between the Markov renewal process and the
multitype branching system.
We introduce the genealogical tree ${\cal T}$ of an individual, which comprises information on the
individual's offspring genealogy, such as family relationships, mutations, death and birth times
of individuals. For $t > 0$, let ${\cal T}_t$ denote the genealogical tree restricted to the time
interval $[0,t]$. Finally, ${\cal T}_t^r$ stands for the reduced tree  obtained from ${\cal T}_t$
by deleting the ancestry lines of those particles, which die before $t$.
We write $\p_i$ for the law of ${\cal T}$, if the process started from
an ancestor of type $i$ with age 0. From the context it will be always clear
when $\p_i$ refers to the branching particle system, or to the Markov
renewal process.

For any given $t > 0$ and ancestry line $w \in {\cal T}_t$, let $t_j(w) \geq 0$
be the total time up to $t$ that $w$ spends in type $j \in \K$.
Introduce the variable
\begin{equation} \label{def_mu}
\mu_j(t) = \min_{w \in {\cal T}_t^r} t_j(w),
\end{equation}
which is the minimal time spent in type $j$ among those particles that are alive at time $t$,
with the usual convention that $\min \emptyset = \infty$. We also define the maximum
spent time in type $j$ up to time $t$:
$$
\sigma_j(t) = \max_{w \in {\cal T}_t} t_j(w).
$$
(Notice that, in this case, the population procreated by the ancestor is not necessarily alive at
time $t$). Let
$$
\nu_{i,j} (t,a) = \p_i \{ \mu_j(t) \leq a \}
$$
denote the probability that starting from $i$, there is a particle at time
$t$, who spent less than $a$ time in $j$.
Note that for $ t < a <  \infty$,
$$
\nu_{i,j}(t, a) = \p_i \{ \mu_j(t) \leq a \}=
\p_i \left \{ \textrm{\,the process is not extinct at } t \right \}
\to 0,
$$
as $t \to \infty$, and for arbitrary $a < \infty$,
$$
\nu_{i,j}(t, a) = \p_i \{ \mu_j(t) \leq a \} \leq
\p_i \left \{ \textrm{\,the process is not extinct at } t \right \}
\to 0,
$$
as $t \to \infty$. Then by a renewal argument we obtain, for $a < t$, that
\begin{eqnarray} \label{rennu}
\nu_{i,i}(t,a)
& = & \int_0^a \Gamma_i (\d s)
\left[ 1 - f_i(\unit - \nu_{\cdot, i}(t-s, a-s) ) \right] \\
\nu_{i,j}(t,a)
& = & 1 - \Gamma_i(t) + \int_0^t \Gamma_i(\d s)
\left[ 1 - f_i( \unit - \nu_{\cdot, j}(t-s, a)) \right]. \nonumber
\end{eqnarray}
Since $1 - f_i(\unit - \z) \leq \sum_{k=1}^K m_{i,k} z_k$, $\z=(z_1,\ldots,z_K)$, we can compare
the solution of \eqref{rennu} with the solution of the linear version
\begin{eqnarray} \label{alpha}
\alpha_{i,i}(t,a)
& = & \int_0^a \Gamma_i(\d s)
\sum_{k=1}^K m_{i,k} \alpha_{k,i}(t-s, a-s) \\
\alpha_{i,j}(t,a)
& = & 1 - \Gamma_i(t) + \int_0^t \Gamma_i(\d s)
\sum_{k=1}^K m_{i,k} \alpha_{k,j}(t-s, a).  \nonumber
\end{eqnarray}
Notice that renewal argument implies again that
$\alpha_{i,j}(t) = \p_i \{ \overline t_j(t) \leq a \}$ is the solution of
the equation system \eqref{alpha}.
Let $\alpha_{i,j}^{(0)}(t,a) = \nu_{i,j}(t,a)$, and let
$\alpha^{(n)}=(\alpha_{i,j}^{(n)})_{i,j=1,\ldots,K}$, where
\begin{eqnarray*}
\alpha_{i,i}^{(n+1)}(t,a)
& = & \int_0^a \Gamma_i(\d s)
\sum_{k=1}^K m_{i,k} \alpha_{k,i}^{(n)}(t-s, a-s) \\
\alpha_{i,j}^{(n+1)}(t,a)
& = & 1 - \Gamma_i(t) + \int_0^t \Gamma_i(\d s)
\sum_{k=1}^K m_{i,k} \alpha_{k,j}^{(n)}(t-s, a).
\end{eqnarray*}
By induction it is clear that $\nu_{i,j}(t,a) \leq \alpha_{i,j}^{(n)}(t,a)$ for all $n$.
We show that the iteration converges to the solution $\alpha_{i,j}(t,a)$, and thus
$\nu_{i,j}(t,a) \leq \alpha_{i,j}(t,a)$.
Let us fix a $t >0$, and introduce the notation
$$
|| x - y ||_t = \sup \{ |x_{i,j}(s, u) - y_{i,j}(s, u)| \, : \
i, j \in\K; \, 0 \leq u < s \leq t \}.
$$
Then we get that, for all $n$,
$$
|| \alpha^{(n+1)} - \alpha^{(n)} ||_t \leq
\max_{i \in \K}
\int_0^t \sum_{k=1}^K m_{i,k} || \alpha^{(n)} - \alpha^{(n-1)} ||_t \Gamma_i(\d s)
= \max_{i \in \K} \Gamma_i(t) || \alpha^{(n)} - \alpha^{(n-1)} ||_t,
$$
where we used that the mean matrix $M$ satisfies $\sum_{k=1}^K m_{i,k} = 1$
for all $i$. The last estimation implies  convergence to the solution of \eqref{alpha}; the proof of
uniqueness of solutions of \eqref{alpha}  follows in the same way. Therefore
we showed that $\nu_{i,j}(t,a) \leq \alpha_{i,j}(t,a)$.
Notice that we only have shown  that any solution of the equation system \eqref{rennu} is
dominated by the \emph{unique} solution of \eqref{alpha}, which does not imply
that \eqref{rennu} has a unique solution. We have proved:

\begin{lemma} \label{comparingRMP}
For every $a \in (0,t)$ we have that
$$
\p_i \left\{ \exists w \in {\cal T}_t^r : \, t_j(w) \leq a \right\}
\leq \p_i \left\{ \bar t_j(t) \leq a \right\},
$$
where the left side is for the branching process, while the right is
for the Markov renewal process.
\end{lemma}

Exactly the same way as in \cite{FV} Lemma 10, we can show a similar
bound.

\begin{lemma} \label{lemma10}
For every $a \in (0,t)$ we have that
$$
\p_i \left\{ \exists w \in {\cal T}_t : t_j(w) \geq a  \right\}
\leq  \p_i \left\{ \overline t_j(t) \geq a \right\},
$$
where  the left side is for the branching process, while the right is
for the Markov renewal process.
\end{lemma}

\section{The Markov renewal process} \label{markovren}

In this section we are going to analyze the auxiliary Markov renewal process.
First consider the discrete Markov chain $X_1, X_2, \ldots$ with transition matrix
$ M$, and let $\pp^*=(p_1^*, \ldots, p_K^*)$ denote its stationary distribution.
We have the following large deviation theorem for Markov chains
(\cite{Gold}, Lemma 2.13.):
For all $\delta > 0$ there exist positive constants $C, c$, such that
\begin{equation} \label{mcldt}
\p_i \left\{ \left| \frac{  t_j(n)}{n} - p_j^* \right| > \delta \right\}
\leq C \e^{-c n},\quad i\in\K,
\end{equation}
where $\p_i$ stands for the probability measure, when the chain starts from
position $i$, and $t_j(n)$ is the number of visits to state $j$ among the
first $n$ steps:
$$
t_j(n) = \# \{ l : X_l = j, l =1,2, \ldots, n \}.
$$
Here and below, several  different constants
arise in the calculations whose precise values are not relevant for our purposes. Therefore, for the
reader's convenience  we chose  not to enumerate these constants. Hence the value of a constant may
vary from line to line. In some proofs we use enumerated constants like $k_1, k_2, \ldots$, whose
values are fixed only in the corresponding proof. Finally, we use some global constants
$c_1, c_2, \ldots$, whose values are the same in the whole paper.

Let  $n_t$ denote the number of renewals up to time $t$.
With these notations we may write
\begin{equation} \label{t_j(t)eq}
\overline t_j(t) = \xi_1^{(j)} + \xi_2^{(j)} + \cdots + \xi_{t_j(n_t)}^{(j)}
+ \eta_j(t) = S_{t_j(n_t)}^{(j)} + \eta_j(t),
\end{equation}
where $\xi_1^{(j)},  \xi_2^{(j)}, \ldots $ are iid random variables with common
distribution function $\Gamma_j$, and
$$
\eta_j(t) = \left\{
\begin{array}{ll}
t - Z_{n_t}, & \textrm{if } X_{n_t}=j,  \\
0, & \textrm{otherwise};
\end{array}
\right.
$$
that is $\eta_j(t)$ is non-zero only for one term, and stands for the
spent lifetime. Here $Z_n$ is the sum of the lifetimes up to the $n^{\textrm{th}}$
renewal, and therefore $Z_n$ is the sum of $n$ independent, but not identically
distributed random variables.

\subsection{A long living particle type}

Let $\gamma\in(0,1]$. Assume that
\begin{eqnarray} \label{assumption-on-gamma_j}
&& 1 - \Gamma_1(x) \sim x^{- \gamma}, \quad \textrm{as } x \to \infty \textrm{ and} \\
&& 1 - \Gamma_j(x) \leq A x^{- \eta_j}, \quad j=2, 3, \ldots, K, \nonumber
\end{eqnarray}
where $A > 0$ and $\eta_j > 1$, $j=2, 3, \ldots, K$.
Put $ \eta= \min\{ \eta_j \, : \, j=2, 3, \ldots, K \}$.
We will show that with high probability the process spends
$c \, t$ times in type $1$.

\begin{lemma} \label{ldtfor1}
There exists $ c_1 > 0$ such that for every $i\in\K$
and $t > 1$
$$
\p_i \left\{ \frac{\overline t_1(t)}{t} \leq c_1 \right\} \leq C t^{1 - \eta},
$$
for some $C > 	0$.
\end{lemma}

\begin{proof*}
For simplicity we omit the lower index $i$.
Recall that $n_t$ stands for the number of renewals up to time $t$. For any $k_2 > 0$
we may write
$$
\p \left\{ \frac{\overline t_1(t)}{t} \leq c_1 \right\}
= \p \left\{ \frac{\overline t_1(t)}{t} \leq c_1 , \, n_t > k_2 \, t \right\} +
\p \left\{ \frac{\overline t_1(t)}{t} \leq c_1 , \, n_t \leq k_2 \, t \right\}.
$$
The first term is easy to estimate.
Due to \eqref{mcldt}, with probability $ \geq 1 - C \, \e^{- c\, t}$ we have
$t_1(n_t) / n_t \geq p_1^* / 2$, and so on this set
$$
\frac{\overline t_1(t)}{t}  \geq \frac{ S^{(1)} _{ t_1(n_t) } }{t_1(n_t)}
\frac{t_1(n_t)}{n_t} \frac{n_t}{t} \geq
\frac{ S^{(1)} _{ t_1(n_t) } }{t_1(n_t)} \frac{k_2 \, p^*_1}{2}.
$$
Truncation and Cram\'er's large deviation theorem shows that
for any $d \in (0, \infty)$ there exist $C, c > 0$, such that for $n \in \N,$
$$
\p \left\{ \frac{ S^{(1)}_n }{n} \leq d \right\} \leq C \, \e^{- c\, n }.
$$
Applying this with $n \sim c t$, the estimation above shows that the
first term $\leq C \, \e^{- c \, t}$ (for some other pair of constants $C, c$)
for any choice of $c_1, k_2$.

Now let us investigate the second term. Clearly $\overline t_1(t) \leq c_1 t $ implies that
$\overline t_j(t) > k_3 t $ for some $j \geq 2$, with $k_3 = (1 - c_1)/ K$. If $n_t \leq k_2 t$
then $\overline t_j(t) \leq S^{(j)}_{t_j(n_t) +1 } \leq S^{(j)}_{ \lfloor k_2 t \rfloor + 1}$
by \eqref{t_j(t)eq},
therefore the probability in question is less then
$$
\p \left\{ \frac{ S^{(j)}_{\lfloor k_2 t \rfloor +1}}{t} > k_3 \right\}
\leq c \, t^{1 - \eta_j} A,
$$
which proves our lemma.
In the last step we used Theorem 2 of Nagaev \cite{nagaev}, which says that for any $c > 0$ and
$x \geq c n$,
\begin{equation} \label{nagaev-ineq}
\p \left\{ S^{(j)}_n - n \mu_j \geq x \right\} \leq 2 n x^{-\eta_j} A
\end{equation}
for $n$ large enough.
In particular, for any $\delta > 0$,
$$
\p \left\{ \frac{ S^{(j)}_n  - n \mu_j}{n} \geq \delta \right\} \leq c n^{1 -\eta_j} A.
$$
\end{proof*}

Notice that by using Lemma \ref{ldtfor_n_t} below we obtain
a stronger result. Namely, for any $ \varepsilon > 0$
$$
\p \left\{ \frac{\overline t_1(t)}{t} \leq c_1 \right\} \leq
t^{\gamma + \varepsilon - \eta},
$$
for $t$ large enough.
Combining this with Lemma \ref{comparingRMP} we obtain
\begin{lemma} \label{ldtfor_occtime1}
For any $ i \in \K$  there is a
$c_1$ and $C > 0$ such that
$$
\p_i \left\{ \frac{\mu_1(t)}{t} \leq c_1 \right\} \leq C t^{1 - \eta}.
$$
\end{lemma}

\subsection{Occupation times for $j \geq 2$}

To analyze the occupation times $\overline t_j(t)$ for $j=2,\ldots,K$, we need
a precise asymptotic for the number of renewals $n_t$.

We start by describing the asymptotic behavior of
$S_n = \xi_1 + \cdots + \xi_n$, the sum of $n$ independent random variables
with distribution function $\Gamma$, for which  $1 - \Gamma(x) \sim x^{- \gamma}$,
$\gamma \in (0,1]$.
In the following, limits of sequences are meant  as $n\to\infty$. We use the same convention for the
continuous parameter $t$.

\begin{lemma} \label{ldtfor_gamma}
Assume that $d_n \to \infty$ if $\gamma < 1$, and
$\log n / d_n \to 0$ when $\gamma = 1$. We have
$$
\p \{ S_n > n^{1/\gamma} d_n \} \leq (1+ o(1)) d_n^{-\gamma}.
$$
Moreover, for $\gamma < 1$
there exist constants $c_\gamma$ such that for any sequences $c_n$ for which
$c_n \to 0$ and $n^{\gamma^{-1} - 1} c_n \to \infty$, the following holds
$$
\p \left\{ S_n \leq c_n n^{1/ \gamma} \right\}
\leq 2 \, \exp \left\{ - \frac{c_n^{- \frac{\gamma}{1 - \gamma}} }
{c_\gamma} \right\}.
$$
If there exists a constant $L > 0$ such that $ \sup_{n} n^{\gamma^{-1} - 1} c_n < L$,
$\gamma \in (0,1]$, then for some $c > 0$
$$
\p \left\{ S_n \leq c_n n^{1/ \gamma} \right\}
\leq \e^{ - c \, n  },
$$
for $n$ large enough.
\end{lemma}

Note that in the case $\gamma = 1$ we can choose $c_n \equiv c > 0$ arbitrary large.
We will use this remark in the proof of Lemma \ref{ldtfor_n_t}.

\begin{proof*}
Let $\{t_n\}$ be a sequence of positive numbers such that $n[1-\Gamma(t_n)]\to0$ if $\gamma<1$ (i.e. $t_n=n^{1/\gamma} d_n$
for some $d_n \to \infty$), and $n \log t_n / t_n \to 0$ when $\gamma=1$ (that is, $t_n= n d_n$, where
$\log n / d_n \to 0$).
Using a theorem of Cline and Hsing (\cite{CH}, Theorem 3.3)
we get that
$$
\lim_{n \to \infty} \sup_{s \geq t_n} \left|
\frac{\p \{ S_n > s \}}{n  [1 - \Gamma(s)] } -1 \right| = 0.
$$
This follows immediately from \cite{CH} if $\gamma < 1$, while for $\gamma =1$
one has to check that the sequence which has to converge to 0, is
$$
\frac{n}{t_n} \int_1^{t_n} x \d \Gamma(x) \sim \frac{n \log t_n}{t_n}.
$$
Writing $s=t_n$ in the form $t_n = n^{1/\gamma} d_n$ we obtain the first statement.

Now we turn to the upper estimates for $S_n/ n^{1/\gamma}$.

We use a truncation method. For $a > 0$ let denote
$$
\xi^{(a)}= \left\{
\begin{array}{ll}
\xi, & \textrm{if } \xi \leq a, \\
a, & \textrm{otherwise},
\end{array} \right.
$$
the truncated variable at $a$. For the first moment of this variable,
as $ a \to \infty$ we have
$$
\mu_a = \E \xi^{(a)}
=  \int_0^a x \d \Gamma(x) + a [ 1 - \Gamma(a) ]
\sim  \left\{
\begin{array}{ll}
\frac{1}{1- \gamma} a^{1- \gamma},
& \textrm{for } \gamma \ne 1, \\
\log a , & \textrm{for } \gamma=1.
\end{array} \right.
$$
For the proof we need Bernstein's inequality  (cf. p. 855 of Shorack and Wellner \cite{SW}).
Let $X_1, X_2, \ldots$ be iid random variables with
$\E X_1 = 0$, and let $\kappa > 0$ and $v > 0$ be constants such that
$\E |X^m| \leq v \kappa^{m-2} m!/2$. Then for the partial sum $S_n=X_1+\cdots+X_n$ the following holds:
\begin{equation} \label{bernstein}
\p \left\{ | S_n | > t \right\} \leq
2 \, \exp \left\{ - \frac{t^2}{2 v n + 2 \kappa t } \right\}.
\end{equation}
Easy computations show that in our case (that is if $ 1 - \Gamma(x) \sim x^{- \gamma} $)
for $ m \geq 2$
$$
\E [ \xi^{(a)} ]^m \sim \frac{m}{m - \gamma} a^{m - \gamma},
\quad \textrm{as } a \to \infty.
$$
Since $ \E | \xi^{(a)} - \mu_a |^m \leq \E [ \xi^{(a)} ]^m + \mu_a^m$,
this shows that in Bernstein's inequality \eqref{bernstein} we
can choose $v = 2 \, a^{2 - \gamma}$ and $\kappa = a$.
Obviously $\p \left\{ S_n \leq c_n n^{1/ \gamma} \right\}
 \leq  \p \left\{ S_n^{(a)} \leq c_n n^{1/\gamma} \right\} $, and so we may write
$$
\p \left\{ S_n \leq c_n n^{1/ \gamma} \right\}
 \leq  \p \left\{ S_n^{(a)} - n \mu_a  \leq c_n n^{1/\gamma} - n \mu_a \right\}
 \leq  \p \left\{ | S_n^{(a)} - n \mu_a | \geq n \mu_a - c_n n^{1/\gamma} \right\},
$$
where in the last inequality we assumed that $ n \mu_a > c_n n^{1/\gamma}$. Applying
Bernstein's inequality with $v= 2 \, a^{2 - \gamma}$
and $\kappa = a$ we obtain
\begin{equation}
\label{X*X}
\p \left\{ S_n \leq c_n n^{1/ \gamma} \right\}
\leq 2 \, \exp \left\{ - \frac{( n \mu_a - c_n n^{1/\gamma})^2}
{4 a^{2 - \gamma} n + 2 a ( n \mu_a - c_n n^{1/\gamma} ) } \right\}.
\end{equation}
Let $n \mu_a = 2 n^{1/ \gamma} c_n$, that is
$
a \sim \left[ 2 ( 1 - \gamma) \right]^{\frac{1}{1-\gamma}}
n^{\frac{1}{ \gamma}} c_n^{\frac{1}{1 - \gamma}}=:a_n.
$
By our assumptions,  $ a_n$ tends to $\infty$. Then
the numerator in the exponential of \eqref{X*X} is $n^{2/ \gamma}c_n^2$, while the
denominator is
$$
4 (\left[ 2 ( 1 - \gamma) \right]^{\frac{1}{1-\gamma}}
n^{\frac{1}{ \gamma}} c_n^{\frac{1}{1 - \gamma}})^{2 - \gamma} n
+ 2 (\left[ 2 ( 1 - \gamma) \right]^{\frac{1}{1-\gamma}}
n^{\frac{1}{ \gamma}} c_n^{\frac{1}{1 - \gamma}})  c_n n^{1/\gamma}
= c_\gamma \, n^{2 / \gamma} c_n^{\frac{2 - \gamma}{1 - \gamma}},
$$
with $ c_\gamma = (10 - 8 \gamma) [2 - 2 \gamma]^{\frac{1}{1 - \gamma}} $.
In this way we obtain finally that
$$
\p \left\{ S_n \leq c_n n^{1/ \gamma} \right\}
\leq 2 \, \exp \left\{ - \frac{c_n^{- \frac{\gamma}{1 - \gamma}} }
{c_\gamma} \right\},
$$
which is the desired bound.
The last assertion in the lemma follows easily from Cram\'er's large deviation
theorem, together with the truncation method.
\end{proof*}

Next we investigate the asymptotic behavior of the number of renewals $n_t$ in our
Markov renewal process, where the lifetime distributions
$\Gamma_1, \ldots, \Gamma_K$ are as in \eqref{assumption-on-gamma_j}.
We show that $n_t$ asymptotically behaves
like the number of renewals $\tilde n_t$ in a \emph{standard} renewal process,
where the tail of the lifetime distribution is $ \sim x^{-\gamma}$.

\begin{lemma} \label{ldtfor_n_t}
Let $c(\cdot)$ be a function such that
\begin{equation}
\label{c(.)}
\lim_{t\to\infty}c(t)=0\mbox{ if $\gamma<1$,}\quad\mbox{and}\quad\lim_{t\to\infty}c(t)\log t=0
\mbox{ if $\gamma=1$}.
\end{equation}
Then for all $i \in \K$,
for $t>0$ large enough,
$$
\p_i \left\{ \frac{n_t}{t^\gamma} \leq c(t) \right\} \leq 2 c(t),
\quad 0 < \gamma \leq 1.
$$
If $\gamma < 1$ then for any $ a \in (0, 1 - \gamma)$, for all $i \in \K$ and all
$t>0$ large enough,
$$
\p_i \left\{ \frac{n_t}{t^\gamma} > t^a \right\}
\leq  2 \, \exp \left\{ - c \, t^{\frac{a}{1 - \gamma} } \right\}.
$$
If $\gamma = 1$, then for any $a > 0$ and $i \in \K$ there is a constant $c > 0$ such
that for all $t>0$ large enough,
$$
\p_i \left\{ \frac{n_t}{t} > a \right\} \leq \e^{- c t}.
$$
\end{lemma}

\begin{proof*}
We drop the lower index $i$.
In order to get an upper bound for $n_t$ let us define
$$
\tilde \Gamma(x) = \prod_{i=1}^K \Gamma_i(x)\,,
$$
which is the distribution function
of the lifetime $\tilde \xi \ed \max\{ \xi^{(1)}, \ldots, \xi^{(K)} \}$,
where $\xi^{(1)}, \ldots, \xi^{(K)}$ are independent and distributed as
$\Gamma_1, \ldots, \Gamma_K$ respectively. Consider a standard renewal process
$\tilde n_t$, $\tilde S_n$ with this lifetime distribution.
Recall that $Z_n$ is the sum of the lifetimes up to the $n^\textrm{th}$ renewal.
Clearly
$$
\p \left \{ n_t \leq t^\gamma c(t) \right\} =
\p \left\{ Z_{\lfloor t^\gamma c(t) \rfloor } > t \right\}
\leq \p \left\{ \tilde S_{\lfloor t^\gamma c(t) \rfloor } > t \right\}.
$$
According to Lemma \ref{maxXY} below, $1 - \tilde \Gamma(x) \sim x^{-\gamma}$. Therefore using
Lemma \ref{ldtfor_gamma} we can write
\begin{eqnarray*}
\p \left\{ \tilde S_{\lfloor t^\gamma c(t) \rfloor } > t \right\}
\leq 2 \, c(t),
\end{eqnarray*}
where we have
assumed  that for $ \gamma < 1$ the convergence
$ t / [ t^\gamma c(t) ]^{1/\gamma} \to \infty $ holds, which is equivalent
to $c(t) \to 0$, and that $c(t) \log t \to 0$ when $\gamma =1$.
This proves the first statement.

To obtain the lower bound we use the simple estimation $Z_n \geq S_{t_1(n)}^{(1)}$,
that is, we simply drop the lifetimes with finite mean. First consider the
case $\gamma < 1$.
Then we have
\begin{eqnarray*}
\left\{  Z_{\lfloor t^{\gamma + a } \rfloor} \leq t \right\}
& \subset &  \left\{  S_{t_1( \lfloor t^{\gamma + a}  \rfloor )}^{(1)} \leq t \right\} \\
& \subset &
\left\{ S_{ \lfloor (p_1^* - \varepsilon) t^{\gamma + a} \rfloor }^{(1)} \leq t, \
\frac{ t_1( \lfloor t^{ \gamma + a }  \rfloor )}{ \lfloor t^{\gamma + a} \rfloor}
> p_1^* - \frac{\varepsilon}{2} \right\}\cup \left\{ \frac{ t_1( \lfloor t^{\gamma + a}  \rfloor )}
{\lfloor t^{ \gamma + a}  \rfloor} \leq p_1^* - \frac{\varepsilon}{2} \right\},
\end{eqnarray*}
hence using Lemma \ref{ldtfor_gamma} and \eqref{mcldt} we have
\begin{eqnarray*}
\p \left\{ \frac{n_t}{t^\gamma} > t^a \right\}
& = & \p \left\{ n_t  > \lfloor t^{\gamma + a}  \rfloor \right\}
= \p \left\{ Z_{\lfloor t^{ \gamma + a} \rfloor} \leq t \right\} \\
& \leq &
\p \left\{ S_{ \lfloor (p_1^* - \varepsilon)  t^{ \gamma + a} \rfloor }^{(1)} \leq t \right\}
+ c \, \e^{- c t^{ \gamma + a} }  \\
& = & \p \left\{ \frac{S_{ \lfloor (p_1^* - \varepsilon) t^{ \gamma + a} \rfloor }^{(1)}}
{ \lfloor (p_1^* - \varepsilon) t^{ \gamma  + a} \rfloor^{1/\gamma} } \leq
\frac{t}{\lfloor (p_1^* - \varepsilon) t^{ \gamma  + a} \rfloor^{1/\gamma}} \right\}
+  c \, \e^{- c t^{\gamma + a} }  \\
& \leq & \p \left\{  \frac{S_{ \lfloor (p_1^* - \varepsilon)  t^{\gamma + a} \rfloor }^{(1)}}
{ \lfloor (p_1^* - \varepsilon)  t^{ \gamma + a } \rfloor ^{1/\gamma} } \leq
c \, t^{-a/\gamma} \right\}
+ c \, \e^{- c t^{ \gamma + a } }  \\
& \leq &
2 \, \exp \left\{ - c \, t^{\frac{a}{1 - \gamma} } \right\}
+ c \, \e^{- c t^{ \gamma + a } }.
\end{eqnarray*}
Taking into account that $\gamma + a > a/(1 - \gamma)$,
we obtain the statement.
Finally, when $\gamma = 1$ we use exactly the same method. Using the
event-decomposition as before, the last part of Lemma \ref{ldtfor_gamma}
and \eqref{mcldt} we have
$$
\p \left\{ \frac{n_t}{t} > a \right\}
 \leq  \p \left\{ S^{1}_{\lfloor a(p_1^* - \varepsilon) t \rfloor }
\leq t \right\} + c \, \e^{- c t }
 \leq  c \, \e^{- c \, t },
$$
thus proving the last assertion of the lemma.
\end{proof*}

The preceding results allow us to obtain the following  estimations for the probabilities of
the smallness and largeness of $\overline t_j(t)/t^\gamma$.

\begin{lemma} \label{ldtfor_type2}
Assume that \eqref{assumption-on-gamma_j} holds, and let
$c(\cdot)$ be a function satisfying \eqref{c(.)}.
Then for all $i\in\K$, any $j = 2,3, \ldots, K$,
every $\varepsilon > 0$ and all $t$ large enough,
$$
\p_i \left\{ \overline t_j(t) \leq t^\gamma c(t) \right\}
\leq  c\, ( c(t) + t^{\varepsilon - \gamma} ), \quad
0 < \gamma \leq 1.
$$
If $\gamma < 1$ then for any $0 < a < 1 -\gamma$,
$$
\p_i \left\{ \overline t_j(t) \geq t^{\gamma + a} \right\}
\leq t^{1 - \eta \gamma},
$$
while for $\gamma = 1$,
$$
\p_i \left\{ \overline t_j(t) \geq a \, t \right\}
\leq c \, t^{1 - \eta}
$$
for any $a > 0$.
\end{lemma}

\begin{proof*}
As before, for simplicity we omit the lower index $i$.
Clearly, for any $\varepsilon > 0$,
$$
\p \left\{ \overline t_j(t) \leq t^\gamma c(t) \right\} \leq
\p \left\{ \overline t_j(t) \leq t^\gamma c(t), n_t \geq t^\varepsilon \right\}
+ \p \left\{ n_t < t^\varepsilon \right\}.
$$
Due to \eqref{mcldt}, on the set $\{ n_t \geq t^\varepsilon \}$  we have, for any $\delta > 0$, that
$t_j(n_t)/ n_t \in (p_j^* - \delta, p_j^* + \delta)$ with probability
$\geq 1 - c\, \e^{- c \, t^\varepsilon}$.
Truncation method and Cram\'er's large deviation theorem show that
$$
\p \left\{ \frac{S^{(j)}_{t_j(n_t)}}{t_j(n_t)} <  \frac{\mu_j}{2},
n_t \geq t^\varepsilon \right\} \leq
c \, \e^{- c \, t^\varepsilon}.
$$
Since
\begin{equation} \label{t_j(t)dec}
\frac{ \overline t_j(t)}{t^\gamma} = \frac{S_{t_j(n_t)}^{(j)} + \eta_j(t)}{t_j(n_t)}
\frac{t_j(n_t)}{n_t} \frac{n_t}{t^\gamma}
\end{equation}
we have obtained that
$$
\p \left\{ \overline t_j(t) \leq t^\gamma c(t), n_t \geq t^\varepsilon \right\}
  \leq    \p \left\{ n_t  \leq c \, t^\gamma c(t)\right\} +
c \, \e^{ - c\, t^\varepsilon}
  \leq   c\, ( c(t) +\e^{ - c\, t^\varepsilon} ),
$$
where in the last step we used Lemma \ref{ldtfor_n_t}. Taking into account that
$\p \left\{ n_t < t^\varepsilon \right\} \leq c \, t^{\varepsilon - \gamma}$,
which follows again from Lemma \ref{ldtfor_n_t}, the first inequality is proved.

For the second part we use a similar technique. We first deal with the case $\gamma < 1$.
For any $0 < b < 1$ we may write
$$
\p \left\{ \overline t_j(t) \geq t^{\gamma + a}  \right\} =
\p \left\{ \overline t_j(t) \geq t^{\gamma + a}, n_t \geq t^b \right\}
+ \p \left\{ \overline t_j(t) \geq t^{\gamma + a},
n_t < t^b \right\}.
$$
If $n_t < t^{b}$,
using that $S_{t_j(n_t)}^{(j)} + \eta_j(t) \leq S_{t_j(n_t)+1}^{(j)}$ we get
$\overline t_j(t) \leq S^{(j)}_{t_j(n_t) +1}
\leq S^{(j)}_{\lfloor t^b \rfloor + 1}$.
Therefore the tail probabilities of the right-hand term satisfy the inequality
\begin{equation} \label{ldtfor_type2_proof}
\p \left\{ S_{ \lfloor t^b \rfloor + 1}^{(j)} \geq
t^{\gamma + a}\right\} \leq
2^{1 + \eta_j} \left( \lfloor t^b \rfloor + 1 \right) t^{- \eta_j (\gamma + a)}
A \leq   t^{1 - \eta \gamma}
\end{equation}
for all $t$ large enough,
where we used again Nagaev's result \eqref{nagaev-ineq}. Note that
we only needed that $\gamma + a > b$.

For the estimation of the other term we use again the decomposition \eqref{t_j(t)dec}.
Since $n_t \geq t^b$, by \eqref{mcldt} we have
$t_j(n_t)/ n_t \in (p_j^*/2, 2 p_j^*)$ with probability
$\geq 1 - c\, \e^{- c \, t^b}$. Due to Lemma \ref{ldtfor_n_t}, for any
$\varepsilon^\prime < a/2$
$$
\p \left\{ \frac{n_t}{t^\gamma} \geq t^{\varepsilon^\prime} \right\}
\leq 2 \, \exp \left\{ - c t^{\frac{\varepsilon^\prime}{1 - \gamma}} \right\}.
$$
Since the orders of these terms are smaller than that of $t^{1-\eta \gamma}$ (given in the statement),
we can work on
$\{ t^b \leq n_t < t^{\gamma + \varepsilon^\prime} \} \cap
\{ t_j(n_t) / n_t \in (p_j^*/2, 2 p_j^*) \}$. On this event, by
\eqref{t_j(t)dec}
$$
\frac{\overline t_j(t)}{t^\gamma} \leq
\frac{S_{t_j(n_t)+1}^{(j)}}{t_j(n_t)}
\frac{t_j(n_t)}{n_t} \frac{n_t}{t^\gamma} \leq
\frac{S_{t_j(n_t)+1}^{(j)}}{t_j(n_t)} 2 p_j^* t^{\varepsilon^\prime},
$$
and so $\overline t_j(t) \geq t^{\gamma + a}$ implies
$S_{t_j(n_t)+1}^{(j)}/t_j(n_t) \geq t^{a - \varepsilon^\prime}/(2 p_j^*).$
Thus for $t$ large enough
\begin{eqnarray*}
&& \vspace{15pt} \p \left\{ \frac{t_j(t)}{t^\gamma} > t^a,
\frac{t_j(n_t)}{ n_t} \in (p_j^*/2, 2 p_j^*),
t^b \leq n_t \leq t^{\gamma + \varepsilon^\prime} \right\} \\
&& \vspace{15pt} \leq \p \left\{ \frac{S^{(j)}_{t_j(n_t)+1}}{t_j(n_t)}
\geq t^{a - \varepsilon^\prime},
t_j(n_t) \geq t^b p_j^*/2 \right\} \\
&& \vspace{15pt} \leq 2 \, t^{1 - \eta ( a + b - \varepsilon^\prime)},
\end{eqnarray*}
where the last inequality follows again from \eqref{nagaev-ineq}.
Choosing $b$ such that $ a + b > \gamma $ and
$\varepsilon^\prime < a + b - \gamma$ we
obtain the desired order $t^{1- \eta \gamma}$.
This, together with \eqref{ldtfor_type2_proof} gives the statement.

The proof in the case $\gamma = 1$ follows a similar approach. For any $b > 0$
$$
\p \left\{ \overline t_j(t) \geq a \, t \right\} \leq
\p \left\{ \overline t_j(t) \geq a \, t, n_t \leq b \, t \right\}
+ \p \left\{ n_t >  b \, t \right\}.
$$
The second summand in the right of the above inequality is exponentially small for any $b > 0$, and
we have already shown in the proof of Lemma \ref{ldtfor1}
that the first one  is less than $c t^{1- \eta}$, provided that $b$ is small enough.
\end{proof*}

Combining the last result with Lemma \ref{comparingRMP} and Lemma \ref{lemma10} respectively
we obtain

\begin{lemma} \label{ldtfor_occtime2}
Assume that \eqref{assumption-on-gamma_j} holds and that
$c(\cdot)$ is a function satisfying \eqref{c(.)}.
Then for all $i\in\K$, all $j = 2,3, \ldots, K$, each
$\varepsilon > 0$ and all $t$ large enough,
$$
\p_i \left\{ \exists w \in {\cal T}^r_t : t_j(w) \leq t^{\gamma} c(t)  \right\}
= \p_i \left \{ \mu_j (t) \leq t^\gamma  c(t) \right\} \leq
c\, ( c(t) + t^{\varepsilon - \gamma} ).
$$
If $\gamma < 1$, then for any $0 < a < 1 - \gamma $
$$
\p_i \left\{ \exists w \in {\cal T}_t : t_j(w) \geq t^{\gamma + a}  \right\} =
\p_i \left \{ \sigma_j (t) \geq t^{\gamma + a }  \right\} \leq
 t^{1 - \eta \gamma },
$$
while if $\gamma = 1$, then for any $a > 0$
$$
\p_i \left\{ \exists w \in {\cal T}_t : t_j(w) \geq a \, t  \right\} =
\p_i \left \{ \sigma_j (t) \geq a \, t  \right\} \leq
c \, t^{1 - \eta}.
$$
\end{lemma}

\section{Extinction results} \label{extinction}

In this final section, we apply the results on
occupation times proved earlier in the paper to analyze extinction properties of our branching particle
system. Let $N_t$ denote the particle system at time $t$, i.e. $N_t$ is
the point measure on $\R^d\times\K$ determined by the positions and types of individuals alive at time
$t\ge 0$. We write $N_t^i$ for the point measure representing the population of type-$i$
particles at time $t$, that is $N_t^i(A) = N_t( A \times \{ i \} )$ for any $A \subset \R^d$,
hence $N_t = N_t^1+\cdots+ N_t^K$.
As before the lower indices in $\p$ and $\E$ refer to the initial
distribution. In particular, $\p_{x,i}$ and $\E_{x,i}$  refer to a population having an ancestor $\delta_{(x,i)}$ of type $i\in\K$, initially at position $x\in \Rd$.

Let $\h: \R^d \times \K \to [0,\infty)$ be continuous function with compact support.  We write
$\langle \mu , \h \rangle = \int \h \, \d \mu$ for any measure  $\mu$ on ${\cal B}(\Rd\times\K)$. Without danger of confusion we also write $\langle\x, \y\rangle  = \sum_{i=1}^K x_i y_i$ for the scalar product
of vectors $\x=(x_1,\ldots,x_K)$ and $\y=(y_1,\ldots,y_K)$.
Assume that the initial population $N_0$ is a Poisson process with intensity measure
$\Lambda = \lambda_1 \, \ell\delta_{\{1\}}+ \cdots + \lambda_K \, \ell\delta_{\{K\}}$,
where $\ell$ is $d$-dimensional Lebesgue measure, and $\lambda_i$, $i\in\K$, are non-negative
constants.

The Laplace transform of our branching process is, for any $t\ge0$, given by
\begin{eqnarray*}
\E\left[\e^{- \langle N_t, \h \rangle }\right]
& = & \exp \left\{ - \sum_{j=1}^K \lambda_j
\int_{\R^d} \E_{x,j} \left[ 1 - \e^{- \langle N_t, \h \rangle} \right] \d x
\right \} \\
& = & \exp \left\{ - \left \langle \Lambda, \unit - \E_{\cdot, \cdot}
\e^{- \langle N_t, \h \rangle }  \right \rangle \right\}.
\end{eqnarray*}
We put $U_i(\h,t,x) = \E_{x,i} \left( 1 - \e^{- \langle N_t, \h \rangle } \right)$.

To prove extinction of $\{N_t$, $t\ge0\}$ it suffices to show that the Laplace transform
of $N_t$ converges to the Laplace transform of the empty population, and for this it is enough
to verify that
$$
\langle  \Lambda , U_{\cdot}(\h, t,\cdot) \rangle \to 0 \quad \textrm{as } t \to \infty,
$$
which is the same as
$$
U_i^+(\h,t) := \int \left[ \E_{x,i} \left( 1 - \e^{- \langle N_t, \h \rangle } \right)
\right] \d x \to 0 \quad \textrm{as }  t \to \infty \mbox{ for all $i\in\K$}.
$$
Let $B\subset\Rd$ be a ball, and assume that $B \times \K \supset \supp \h$. Then
$$
1 - \e^{- \langle N_t, \h \rangle } \leq I( N_t(B \times \K)  > 0 ),
$$
which implies
$$
\E_{x,i} \left[ 1 - \e^{- \langle N_t, \h \rangle }  \right]
\leq \p_{i,x} \left\{  N_t(B \times \K)  > 0 \right\}.
$$
Conversely, if $\h |_{B \times \K} \geq 1$, then
$$
1 - \e^{- \langle N_t, \h \rangle } \geq (1 - \e^{-1} )
I( N_t(B \times \K)  > 0 ),
$$
and so
$$
\E_{x,i} \left[ 1 - \e^{- \langle N_t, \h \rangle }  \right]
\geq  (1 - \e^{-1} ) \p_{i,x} \left\{  N_t(B \times \K)  > 0 \right\}.
$$
In this way we get that
\begin{lemma} \label{extinct} Extinction of $\{N_t$, $t\ge0\}$ occurs if, and only if for any
bounded Borel set $B\subset\Rd$,
$$
\int_{\R^d}  \p_{i,x} \left\{  N_t(B \times \K)  > 0 \right\} \d x \to 0
\ \textrm{for all } i\in\K, \textrm{ as } t \to \infty.
$$
\end{lemma}

Put $\alpha = \min\{ \alpha_i : i= 1,2, \ldots, K \}$. Recall the following result from  \cite{FV}:

\begin{lemma} (Fleischmann \& Vatutin). \label{L13FV}
For each bounded $B \subset \R^d$
$$
\sup_{t \geq 1} \int_{\R^d  \backslash C(t,L)} \E_{x,i} N_t(B \times \K) \d x
\longrightarrow 0 \textrm{ as } \ L \uparrow \infty,
$$
where $C(t,L) = \{ x \in \R^d : |x| \leq L t^{1/\alpha} \}$.
\end{lemma}

This means that extinction of $\{N_t$, $t\ge0\}$ occurs if, and only if for any
bounded Borel set $B\subset\Rd$, and for $L$ large enough
\begin{equation} \label{exteq}
\int_{C(t,L)}  \p_{i,x} \left\{  N_t(B \times \K)  > 0 \right\} \d x \to 0
\ \mbox{ for all $i\in\K$} \textrm{ as } \ t \to \infty.
\end{equation}

\subsection{Lifetimes with finite means}

When the lifetimes have finite mean and the dimension is small, it is not necessary to analyze the
occupation times in order to prove local extinction. As we are going to show,
in this case a simple estimation and the asymptotics  of the extinction probabilities of critical
multitype branching processes give the result.

Let $F^{(i)}$ denote the probability
generating function of the process starting from a single particle of type $i$:
\begin{equation} \label{defF}
F^{(i)}(t; s_1, \ldots, s_K) =
\E_i\left[ s_1^{N^1_t(\Rd)}\cdots s_K^{N^K_t(\Rd)} \right],\quad 0\le s_j\le 1,\ j \in \K.
\end{equation}
Put $Q^{(i)}(t; s_1, \ldots, s_K) = 1 - F^{(i)}(t; s_1, \ldots, s_K)$ and
$q^{(i)}(t; s) =Q^{(i)}(t; s, \ldots, s)$.
Clearly
$$
\p_{i,x} \left\{  N_t(B \times \K)  > 0 \right\}
\leq \p_{i} \left\{ \textrm{the process is not extinct at time $t$} \right\}.
$$
Consider the discrete--time
multitype Galton--Watson process $\{ \X_n \}$, with the same
offspring distributions as in the branching particle system.
Let  $\v$ and $\u$ respectively denote the  left and right
normed eigenvectors of the mean matrix $M$, which are determined by:
\begin{equation} \label{eigen-def}
\v M = \v,  \quad M \u = \u, \quad \v \u=1, \quad \unit \u =1.
\end{equation}
Since by assumption $M$ is stochastic, $\u = K^{-1} \unit$.
Let $\f_n=(f^1_n, \ldots, f^K_n)$ denote the generating function of the $n^\textrm{th}$
generation, that is $f^i_n(\x) = \E_i\left[ \x^{\X_n}\right]$ and put $\f_1(\x) = \f(\x)$.
It is well-known that $\f_{n+1}(x) = \f (\f_n( \x))$.
Let us assume that
\begin{equation} \label{branchingassumption}
x - \langle \v, \unit - \f(\unit - \u \, x) \rangle \sim x^{1+\beta} L(x) \quad
\textrm{as } x \to 0,
\end{equation}
where
$\beta \in (0,1]$ and $L$ is slowly varying at 0 in the sense that
$\lim_{x \to 0}L(\lambda x)/ L(x) = 1$
for every $\lambda > 0$.
In this case, for the survival probabilities it is known that
$$
\unit - \f_n(0) = ( \u + o(1) ) n^{-1/\beta} L_1(n) \mbox{ as } n\to\infty,
$$
where $L_1$ is slowly varying at $\infty$ (see Theorem 1 in \cite{Vat77} or
Theorem 1 in \cite{Vat78}).
Moreover, assume that
\begin{equation} \label{lifetimeassumption}
\lim_{n \to \infty} \frac{n [1 - \Gamma_i(n)]}{\langle \v, \unit - \f_n(0) \rangle} =0,
\quad i=1,2, \ldots, K.
\end{equation}
Then
\begin{equation} \label{extest}
Q^{(i)}(t; 0) = \p_i \left\{ \textrm{the process is not extinct at $t$}\right\} \sim
u_i \, D^{\frac{1}{\beta}} t^{-\frac{1}{\beta}} L_1(t)\mbox{ as $t \to \infty$},
\end{equation}
where $D= \sum_{i=1}^K u_i v_i \mu_i$; see Theorem 2 in \cite{Vat78}.

Using the estimate above, we obtain the following theorem.

\begin{theorem} \label{finmean}
Assume that \eqref{branchingassumption} and \eqref{lifetimeassumption} hold. Then
for $d < \alpha/ \beta$ the process $\{N_t$, $t\ge0\}$ suffers local extinction.
\end{theorem}

\begin{proof*}
Due to \eqref{extest}, for any $\varepsilon > 0 $
$$
Q^{(i)}(t; 0) \leq c \, t^{-\frac{1 - \varepsilon}{\beta}},\quad i\in\K.
$$
Plugging this into \eqref{exteq} we get
$$
\int_{C(t,L)}  \p_{i,x} \left\{  N_t(B \times \K)  > 0 \right\} \d x \leq
c \, t^{\frac{d}{\alpha} - \frac{1 - \varepsilon}{\beta}}.
$$
Since by assumption $d < \alpha/\beta$,   for some $\varepsilon>0$ the exponent of $t$ in the
above inequality is negative, which implies that the integral in the left-hand side tends to 0.
\end{proof*}\bigskip

\noindent \textbf{Remark 1.} When the generating functions $f_i$, $i\in\K$, are of the form
$
f_i(s, \ldots, s) = f_i(s) = s + c\, (1 - s)^{1 + \beta_i}
$ where $\beta_i\in(0,1]$,
it is easy to verify that \eqref{branchingassumption} holds with
$\displaystyle\beta= \min\{\beta_i:i\in\K\}$, and that \eqref{lifetimeassumption} is fulfilled if for some $\varepsilon > 0$
$$
\lim_{n \to \infty} n^{1 + \frac{1}{\beta} + \varepsilon} [1 - \Gamma_i(n)] =0,
\quad i=1,2, \ldots, K.
$$

\medskip

\noindent \textbf{Remark 2.}
We remark  that  we do not need the precise asymptotic decay of the non-extinction probabilities
given in \eqref{extest}; it suffices to know an asymptotic order of decay. In order to get this,
instead of assuming in \eqref{branchingassumption} that $L(\cdot)$ is slowly varying at 0, it is
enough to suppose that $L$ is an \emph{S--O varying function\/}, meaning that
there exists an $A > 0$ such that $\limsup_{x \to 0} L(\lambda x ) / L(x) < A$ for
any $\lambda > 0$. S--O varying functions were introduced by Drasin and Seneta \cite{DS}.
The definition immediately implies that $\liminf_{ x \to 0} L(\lambda x)/ L(x) > A^{-1}$ for all
$\lambda>0$. It was shown in \cite{DS} that every S--O varying function admits a representation
as the product of a
slowly varying function and a bounded (away from 0 and $\infty$) function.
A careful analysis of the proof of Theorem 1 in \cite{Vat77} shows that, under the
S--O varying assumption on $L$, we have that for any $\varepsilon > 0$ and for all $n$
large enough,
$$
| \unit - \f_n(0)| \leq   n^{- \frac{1 - \varepsilon}{\beta}}.
$$
Since this estimate is precisely what we use in the proofs of our extinction  theorems, all these
results (including the infinite mean case) remain true in this more general setup.
If \eqref{lifetimeassumption} holds (which in particular implies that the
lifetimes have finite mean), we obtain that for any
$\varepsilon > 0$ and for all $t$ large enough,
$$
Q^{(i)}(t; 0) = \p_i \left\{ \textrm{the process is not extinct at $t$}\right\}
\leq t^{- \frac{1 - \varepsilon}{\beta}}.
$$

\subsection{A lifetime with infinite mean -- Case A} \label{infA}

From now on we assume that there is exactly one lifetime distribution
with infinite mean; more precisely we assume \eqref{assumption-on-gamma_j}.
Moreover, in this subsection we additionally assume that
$\alpha = \min \{ \alpha_i : i\in\K \} = \alpha_1$,
that is, the long-living particle type is the most mobile as well.

In the following, $\p^{\theta}_{x,i}$ denotes the distribution of the population starting with a
single individual $\delta_{(x,i)}$  of age $\theta \geq 0$.

\begin{lemma} \label{L1VW} For all $(x,i)\in\Rd\times\K$, all bounded Borel $B\subset\Rd$ and all $t>0$,
$$
\p_{x,i}^\theta \left\{ N_t(B \times \K) > 0 \right\} \leq c_2 \,
\left( t^{-d/\alpha} + t^{1 - \eta} \right),
$$
where the constant $c_2$ is independent of $\theta, x$ and $i$.
\end{lemma}

\begin{proof*}
Put $A= \{ \mu_1(t) \leq c_1 t \}$, i.e.~$A$ is the event that
there exists a branch $ w \in {\cal T}_t^r$
such that $t_1(w) \leq c_1 \, t$, so the process spends less than $c_1 \, t$
time in type 1 for some branch $w$.
Clearly, due to Lemma \ref{ldtfor_occtime1}, we may write
\begin{eqnarray*}
\p \{ N_t(B \times \K ) > 0 \}
& \leq & \p \{ A \} + \p \{ N_t(B \times \K ) > 0, A^c \} \\
& \leq & c t^{1 - \eta} + \E \left[  N_t(B \times \K ) I_{A^c} \right]. 
\end{eqnarray*}
Conditioning on the reduced tree and noting that $A$ is ${\cal T}_t^r$ measurable
we have
\begin{eqnarray*}
\E_{x,i}^\theta \left[  N_t(B \times \K ) I_{A^c} \right] 
& = & \E_{x,i}^\theta \sum_{j=1}^K \sum_{l=1}^{N_t^j(\Rd)} I_{A^c} I( W_j^l(t) \in B ) \\
& = & \E_{x,i}^\theta  \E_{x,i}^\theta \left[
\sum_{j=1}^K \sum_{l=1}^{N_t^j(\Rd)}  I_{A^c} I( W_j^l(t) \in B ) \bigg|
{\cal T}_t^r \right] \\
& = & \E_{x,i}^\theta \left[ I_{A^c} \sum_{j=1}^K \sum_{l=1}^{N_t^j(\Rd)}
\p_{x,i}^\theta  \left\{  W_j^l(t) \in B  \big| {\cal T}_t^r \right\} \right],
\end{eqnarray*}
where, given ${\cal T}_t^r$,
\begin{equation}
\label{sojourn_times}
W_j^l(t) \stackrel{{\cal D}}{=} W(t_1, \alpha_1) + \cdots + W(t_K, \alpha_K).
\end{equation}
Here $t_j$ is the time that a branch of the reduced tree spent in type $j$,
$t_1 + \cdots + t_K = t$,
and $\{W(t, \alpha_j)$, $t\ge0\}$ are independent symmetric $\alpha_j$-stable motions starting from 0, $j=1,\ldots, K$.
Since on the complement of $A$ any branch spent at least $c_1 t$ time in type 1, we have
\begin{eqnarray*}
\p_{x,i}^{\theta} \left\{  W_j^l(t) \in B ; \, A^c  \big| {\cal T}_t^r \right\}
& = & \int p_{t-t_1}(x, \d y) \int_{B - y} p_{t_1}^{\alpha_1}(y, \d z) \\
& \leq & c t_1^{- d / \alpha_1}  = c t^{-d/\alpha}
\end{eqnarray*}
(where $p_{t-t_1}(x,\d y)$ stands for $p^{(\alpha_2)}_{t_2}* \cdots *p^{(\alpha_K)}_{t_K}(x,\d y)$),
and we may continue writing the long equality as
$$
\leq c t^{-d/\alpha} \sum_{j=1}^K  \E N_t^{(j)}(\R^d) \leq  c t^{-d/\alpha}.
$$
Summarizing we obtain
$$
\p_{x,i}^\theta \left\{ N_t(B \times \K) > 0 \right\} \leq c \, t^{-d/\alpha}
+ c \, t^{1 - \eta}.
$$
\end{proof*}

Besides Lemma \ref{L1VW}, our other  key tool is an analogue of Lemma 3 in \cite{VW}. Recall the notations after \eqref{defF}. The proof is an easy multidimensional extension of the proof
in \cite{VW}.

\begin{lemma} \label{L3VW}
If $\eta - 1 > d/ \alpha$, then for any $x \in \R^d, t > 0, i \in \K$ and
$ u \in (0,t- c_2^{\alpha/d})$,
$$
\p_{x,i} \left\{ N_t( B \times \K) > 0 \right\} \leq
q^{(i)}\left(u; 1 - c_2 \, (t - u)^{-d/\alpha}\right),
$$
where the constant $c_2$ is given in Lemma \ref{L1VW}.
\end{lemma}

\begin{proof*} Let $|N_r|\equiv
(N_r^1(\Rd),\ldots,N_r^K(\Rd))$, $r\ge0$.
For any $u \in (0,t)$,
\begin{eqnarray}\label{eesstt}
\p_{x,i} \left\{ N_t( B \times \K) > 0 \right \}
& = & \sum_{\k \neq(0,\ldots,0)
} \p_{x,i} \left\{ |N_u| = \k, N_t(B \times \K) > 0 \right \} \\ \nonumber
& = &  \p_{x,i} \left\{ |N_u| \neq 0 \right\} - \sum_{\k \neq(0,\ldots,0)
} \p_{x,i} \left\{ |N_u| = \k, N_t(B \times \K) = 0 \right \},
\end{eqnarray}
where
$$
 \p_{x,i} \left\{ |N_u| = \k, N_t( B \times \K) = 0 \right \}
  =
\E \left[ \p \left\{  N_t( B \times \K) = 0 \big|  |N_u| = \k, \Theta_\k, Y_\k
\right\} \, I(|N_u| = \k) \right].
$$
Here  $\Theta_\k$ is the vector of ages, and $Y_\k$  the vector of positions of individuals alive
at time $u$. Using independence and Lemma \ref{L1VW},
the conditional probability inside the above expectation gives
\begin{eqnarray*}
&& \p \left\{  N_t( B \times \K) = 0 \big|  |N_u| = \k, \Theta_\k, Y_\k \right\} \\
&& \hspace{60pt} = \prod_{j=1}^K \prod_{l=1}^{k_j}
\p_{y_{l, j}, j}^{\theta_{l,j}}
\left\{ N_{t - u }(B \times \K) = 0 \right\} \\
&& \hspace{60pt} =  \prod_{j=1}^K \prod_{l=1}^{k_j}
\left( 1 - \p_{y_{l, j}, j}^{\theta_{l,j}}
\left\{ N_{t - u }(B \times \K) > 0  \right\} \right) \\
&& \hspace{60pt} \geq  \left( 1 - c_2 \, (t-u)^{-d/\alpha} \right)^{|\k|},
\end{eqnarray*}
where in the last inequality we used that $ 1 - c_2 \, (t-u)^{-d/\alpha} > 0$.
Therefore we obtain
$$
\p_{x,i} \left\{ |N_u| = \k, N_t(B \times \K) = 0 \right \}
\geq \left( 1 - c_2 \, (t-u)^{-d/\alpha} \right)^{|\k|} \p_{i} \left\{ |N_u| = \k \right\}.
$$
Substituting this estimate back into \eqref{eesstt}, we finally get
\begin{eqnarray*}
\p_{x,i} \left\{ N_t(B \times \K) > 0 \right \}
& \leq & \sum_{\k \in \N^d \backslash 0}
\left[ 1 - \left( 1 - c_2 \, (t-u)^{-d/\alpha} \right)^{|\k|} \right] \p_{i}
\left\{ |N_u| = \k \right\} \\
& = & Q^{(i)}(u; 1 - c_2 \, (t-u)^{-d/\alpha}, \ldots, 1 - c_2 \, (t-u)^{-d/\alpha})\\
& = & q^{(i)} (u; 1 - c_2 \, (t-u)^{-d/\alpha}).
\end{eqnarray*}
\end{proof*}

Let us define the set
\begin{equation} \label{Lambda}
\Lambda= \{ \s\in[0,1]^K : \f(\s) \geq \s \},
\end{equation}
where an inequality of the form $(x_1,\ldots,x_K) \geq (y_1,\ldots,y_K)$ means here that
$x_i \geq y_i$ for $i=1,2, \ldots, K$.

We remark that, since
$ \unit - \f(\unit - \u x) \leq M \u x = \u x$, we have $\unit - \u x \in \Lambda$
for all $x$ with $ 0 < \u x \leq \unit$. In our case $\u = K^{-1} \unit$, and this implies that
the diagonal $\{(s,\ldots,s) : s\in[0,1]\}$ is contained in $\Lambda$.

For given matrix families  $A(t) = ( a_{ij}(t) )_{i,j}$ and
$B(t) = ( b_{ij}(t) )_{i,j}$, $t\ge0$, let us define the matrix convolution $C = A * B$ by
$$
c_{ij}(t) = \sum_{k=1}^K \int_0^t a_{ik}(t-s) b_{kj}(\d s).
$$
The convolution of a matrix and a vector is defined analogously. Put
$M^1_\Gamma (t) = ( m_{ij} \Gamma_i(t) )_{i,j}$ and recursively define
$$
M^{n+1}_\Gamma (t) = M^1_\Gamma (t) * M^n_\Gamma (t),\quad n=1,2,\ldots.
$$
Put also  $M^0_\Gamma (t) = ( \delta_{ij} \Gamma_i^0(t) )_{i,j}$, where
$\Gamma_i^0(t)$ is the distribution function of a constant 0 random variable. Notice that
$M^0_\Gamma (t)$ constitutes the unit element in matrix convolution. The following multidimensional
comparison lemma is borrowed from \cite{Vat78}, which is a  generalisation of
Goldstein's comparison lemma \cite{Gold}.

\begin{lemma} \label{comparison}
For any $t > 0$, any natural $n$ and for all $\s \in \Lambda$,
\begin{eqnarray*}
\unit  - \f_n(\s) - M_\Gamma^n * [ (\unit - \s) \otimes \Gamma] (t)
& \leq &
\unit - F(t; \s) \\
& \leq &
\unit - \f_n(\s) + \sum_{j=0}^{n-1} M_\Gamma^j * [ (\unit  - \s) \otimes [\unit - \Gamma]](t).
\end{eqnarray*}
\end{lemma}
Here  $\x \otimes \y := (x_1 y_1, x_2 y_2, \ldots, x_K y_K)$ if $\x=(x_1,\ldots,x_K)$ and $\y=(y_1,\ldots,y_K)$.

We are going to use below the upper bound given in Lemma \ref{comparison}. The following lemma is
Lemma 5 in \cite{Vat79}.

\begin{lemma} \label{twobp}
Consider two critical multitype branching processes sharing the same branching mechanism,
with corresponding lifetime distributions $\Gamma(t) = (\Gamma_1(t), \ldots, \Gamma_K(t))$ and
$\Gamma^*(t) = (\Gamma^*_1(t), \ldots, \Gamma^*_K(t))$. Assume that
$\Gamma(t) \geq \Gamma^*(t)$ for all $t \geq 0$. Then for all $t \geq 0$ and $\s \in \Lambda$,
$$
F(t; \s) \leq F^*(t; \s),
$$
where $F$ and $F^*$ are, respectively,  the vector generating  functions of the number of particles at
time $t$ in the first and second process.
\end{lemma}

The main result in this section is the following theorem.

\begin{theorem} \label{infmeanA}
Assume that \eqref{branchingassumption} holds, the mean matrix $M$ is stochastic, and
the lifetimes satisfy
$1 - \Gamma_1(t) \sim t^{-\gamma}$ for some  constant $\gamma \leq 1$, and
$$
1 - \Gamma_j(x) \leq A \, x^{- \eta_j}, \quad j=2, 3, \ldots, K,
$$
where $\eta_j > 1$, $j=2, 3, \ldots, K$.
Put $ \eta= \min\{ \eta_j \, : \, j=2, 3, \ldots, K \}$. If $ \eta - 1 > d/ \alpha$ and $d < \frac{\alpha \gamma}{ \beta }$,
then the process suffers local extinction.
\end{theorem}

\begin{proof*}
Define the distribution function
$$
\tilde \Gamma (t) = \prod_{i=1}^K \Gamma_i(t),
$$
which is the distribution function of
$\tilde \xi = \max\{ \xi_1, \ldots, \xi_K\}$,
where the random variables $\xi_i$, $i=1,\ldots,K$, are independent with distribution function $\Gamma_i$. Lemma \ref{maxXY}
below shows that $1 - \tilde \Gamma(t) \sim t^{-\gamma}$.
Consider a new branching process
where the branching mechanism is unchanged, but the lifetimes of all types have distribution
$\tilde \Gamma$, and let $\tilde F(t; \s)$ denote its generating function at time $t$.
Clearly, the choice of $\tilde \Gamma$ shows that Lemma \ref{twobp} is applicable, and so
for $\s \in \Lambda$,
\begin{equation} \label{Fineq}
\tilde F(t; \s) \leq F(t; \s).
\end{equation}
(Notice that $\Lambda$, as defined in \eqref{Lambda}, depends only on the branching mechanism
of our process). Now we apply the comparison
lemma for this new process. Since now all the lifetimes have the
same distribution,
$$
M_{\tilde \Gamma}^n (t) = M^n \, \tilde \Gamma^{*n}(t),
$$
where $^{*n}$ stands for the $n$-fold convolution. Moreover, for $\s = s \, \unit$,
$$
M_{\tilde \Gamma}^j * [ (\unit - \s) (\unit - \tilde \Gamma)] (t)
= (1 - s) ( \tilde \Gamma^{*j}(t) - \tilde \Gamma^{*(j+1)}(t) ) M^j \unit =
(1 - s) ( \tilde \Gamma^{*j}(t) - \tilde \Gamma^{*(j+1)}(t) )  \unit,
$$
where we used the simple fact that $M^j$ is stochastic if $M$ is stochastic. Thus, in the rightmost
inequality of Lemma \ref{comparison} we get a telescopic sum, and therefore we obtain
$$
\unit - \tilde F(t; \unit s) = \tilde Q (t; \unit s) \leq
\unit - \f_n(\unit s) + (1 - s) [ 1 - \tilde \Gamma^{*n}(t) ] \unit.
$$
According to \eqref{branchingassumption}, for the survival probabilities we have
$$
1 - f_n^{(i)}(\unit s) \leq 1 - f_n^{(i)}(0) \leq c \, n^{-\frac{1}{\beta}}.
$$
Taking into account \eqref{Fineq} we have, for $ s \in (0,1)$,
$$
Q^{(i)}(t; \unit s) = 1 - F^{(i)} (t; \unit s)  \leq
c \, n^{-\frac{1}{\beta}} + ( 1 - s) \p \left\{ S_n > t \right \}.
$$
Choosing $ n = t^{\gamma/(1+ \varepsilon)}$ and using that, by Lemma \ref{ldtfor_gamma},
$$
\p \left\{ S_n > t \right \} =
\p \left\{ S_n > n^{\frac{1+ \varepsilon}{\gamma}} \right \} \leq
2 \, n^{-\varepsilon} = 2 \, t^{- \gamma \varepsilon/(1+ \varepsilon)},
$$
we get
$$
q^{(i)}(t; 1 -s) = Q^{(i)} (t; (1 - s) \unit ) \leq
c \, t^{-\frac{\gamma}{(1+\varepsilon) \beta}}
+ s \,  t^{- \gamma \varepsilon/(1+ \varepsilon)}.
$$
Hence, choosing $u = t/2$ in Lemma \ref{L3VW} we obtain the inequality
\begin{equation} \label{q-ineq}
q^{(i)}(u; 1 - c_2 \, (t - u)^{-d/\alpha}) \leq
c \, t^{-\frac{\gamma}{(1+\varepsilon) \beta}} +
c\, t^{-\frac{d}{\alpha}}
t^{-\frac{\gamma \, \varepsilon}{1 + \varepsilon}}.
\end{equation}
Multiplying by $t^{d/\alpha}$, the second term in the right of  \eqref{q-ineq} goes to 0, while in
the first one the exponent of $t$ becomes
$$
\frac{d}{\alpha}  -\frac{\gamma}{(1+\varepsilon) \beta},
$$
and this is negative if $d < \alpha \gamma / \beta$ and $ \varepsilon$ is small enough.
\end{proof*}

The simple lemma we used above is the following:

\begin{lemma} \label{maxXY}
Let $X, Y$ be independent non-negative random variables with corresponding distribution functions
$F$ and $G$. Assume that $1 - F(x) \sim x^{-\gamma}$ and $\E Y < \infty$.
Then for the distribution of $Z = \max \{ X, Y \}$ we have
$$
1 - H(z) := \p \left\{ Z > z \right\} \sim z^{-\gamma},
$$
as $ z \to \infty$.
\end{lemma}

\begin{proof*}
Since $\E Y < \infty$, we have $ y (1 - G(y)) \to 0$. Hence,
\begin{eqnarray*}
z^\gamma [1 - H(z)]
& = & z^{\gamma} \left[ 1 - \p \left\{ \max \{ X, Y \} \leq z \right\} \right] \\
& = & z^\gamma \left[ 1 - F(z) G(z) \right]
= z^\gamma \left[ 1 - F(z) + F(z) ( 1 -  G(z) ) \right] \\
& = & z^\gamma \left[ 1 - F(z) \right] +
z^\gamma F(z) \left[ 1 - G(z) \right]  \to 1.
\end{eqnarray*}
\end{proof*}

\subsection{A lifetime with infinite mean -- Case B} \label{infB}

Now let us investigate the case when $\alpha_1$ is not the minimal
$\alpha = \min \{ \alpha_i : i= 1, 2, \ldots, K \}$. Without loss of generality, let us
assume that $ \alpha = \alpha_2$.

Notice that Lemma \ref{L1VW} is true in this case with exponent $- d/ \alpha_1$, and so
the variation of Lemma \ref{L3VW} also remains true. We state it for the easier reference.

\begin{lemma} \label{L3VWwith-nu}
If $\eta - 1 > d/ \alpha_1$, then for any $x \in \R^d, t > 0, i \in \K$ and
$u \in (0,t- c_2^{\alpha_1/d})$,
$$
\p_{x,i} \left\{ N_t( B \times \K) > 0 \right\} \leq
q^{(i)}\left(u; 1 - c_2 \, (t - u)^{-d/\alpha_1}\right),
$$
where the constant $c_2$ is given in Lemma \ref{L1VW}.
\end{lemma}

Put
\begin{equation} \label{effmob}
v = \max \left\{ \frac{1}{\alpha_1}, \frac{\gamma}{\alpha} \right\}.
\end{equation}

\begin{lemma} \label{L13FVwith-v}
Assume that $\gamma \eta > d/\alpha + 1$.
If $\gamma < 1$, then
for any $\varepsilon > 0$, for any $i \in \{ 1,2, \ldots, K\}$ and for
any bounded Borel set $B$,
$$
\lim_{t \to \infty} \int_{|x| \geq t^{v + \varepsilon}}
\p_{x,i} \left\{ N_t(B \times \K) > 0 \right\} \d x = 0.
$$
For $\gamma = 1$ (then necessarily $v = 1/\alpha$),
$$
\lim_{L \to \infty} \limsup_{t \to \infty}
\int_{|x| \geq L t^{v}}
\p_{x,i} \left\{ N_t(B \times \K) > 0 \right\} \d x = 0.
$$
\end{lemma}

\begin{proof*}
Without loss of generality, we will assume that $B$ is a ball with radius $r$ centered at the
origin. First, consider the case $\gamma < 1$.
Put $C(t) = \{ |x | \leq t^{v + \varepsilon} \}$ and let $\varepsilon^\prime < \alpha \varepsilon$.
Recall the definition of
$\sigma_j(t)$ after \eqref{def_mu} and put
$$
A= \{ \sigma_2(t) \leq t^{\gamma + \varepsilon^\prime},
\sigma_3(t) \leq t^{\gamma + \varepsilon^\prime},  \ldots,
\sigma_K(t) \leq t^{\gamma + \varepsilon^\prime }\},
$$
namely $A$ is the set where, for all ancestry lines, the spent time in
type $j$ up to $t$ is less than  $t^{\gamma + \varepsilon^\prime}$ for all $j=2,3, \ldots, K$.

First we work on the set $A^c$. By Lemma \ref{ldtfor_occtime2},
$$
\p \{ A^c \} \leq \sum_{j=2}^K
\p \left\{ \sigma_j(t) > t^{\gamma + \varepsilon^\prime} \right\}
\leq K \, t^{1 -\eta \gamma }.
$$
According to Lemma \ref{L13FV},
$$
\sup_{t \geq 1} \int_{| x | \geq L \, t^{1 / \alpha}}
\E_{x,i} \left[ I(A^c) N_t(B \times \K) \right] \d x \to 0 \quad
\textrm{as } L \to \infty,
$$
hence, it suffices to integrate on the region $t^{v+ \varepsilon} \leq |x| \leq L t^{1/\alpha}$.
On the other hand
$$
\p_{x,i} \left\{ A^c ,  N_t(B \times \K) > 0 \right\} \leq \p \{ A^c \},
$$
and so
$$
\int_{ L t^{1 / \alpha} \geq |x| \geq t^{v + \varepsilon}}
\p_{x,i} \left\{ A^c , N_t(B \times \K) > 0 \right\} \d x \leq
c \, t^{d/ \alpha} t^{ 1 - \eta \gamma } \to 0
$$
due to our assumption.

From now on we work on $A$. Translation invariance of the motion shows that
$$
\int_{\R^d \backslash C(t)} \E_{x,i} I(A) N_t( B \times \K) \,\d x
= \int_{\R^d \backslash C(t)} \E_{0,i} I(A) N_t( (B - x) \times \K) \,\d x.
$$
By conditioning on the reduced tree, we can write
$$
\E_{0,i} I(A) N_t( (B - x) \times \K)
= \E_{0,i} I(N_t  \ne 0) I(A)  \sum_{j=1}^K \sum_{m=1}^{N_t^j(\R^d)}
\p_{0,i} \left\{ W_j^m(t) \in B - x | {\cal T}_t^r \right\},
$$
where $
W_j^m(t) \stackrel{{\cal D}}{=} W(t_1, \alpha_1) + \cdots + W(t_K, \alpha_K)
$ as in \eqref{sojourn_times}.
Integrating we obtain
\begin{eqnarray*}
\lefteqn{\int_{\R^d \backslash C(t)} \E_{0,i} I(A) N_t( (B - x) \times \K) \,\d x}\\ &=&
\E_{0,i} I(N_t  \ne 0) I(A) \sum_{j=1}^K \sum_{m=1}^{N_t^j(\R^d)}
\int_{\R^d \backslash C(t)} \d x \int_{B - x}
\p_{0,i} \left\{ W_j^m(t) \in \d y | {\cal T}_t^r \right\}.
\end{eqnarray*}
Since $|x + y | \leq r$ and $|x| \geq t^{v + \varepsilon}$, we have
$|y| > t^{v + \varepsilon} - r \geq  t^{v + \varepsilon}/2$ for all
$t$ large enough. Using Fubini's theorem and that
$\int_{|x+y | \leq r} \d x =: c(r)$ independently of $y$,  the double
integral can be bounded from above by
$$
c(r) \p_{0,i} \left\{ |W_j^m(t) | \geq \frac{t^{v + \varepsilon}}{2}
\Big| {\cal T}_t^r \right\}.
$$
On the event $A$ we can write
\begin{eqnarray*}
\p_{0,i} \left\{ |W_j^m(t) | \geq \frac{t^{v + \varepsilon}}{2}
\, \Big| \, {\cal T}_t^r \right\}
& \leq & \sum_{k=1}^K \p \left\{ |W(t_k, \alpha_k) | \geq
\frac{t^{v + \varepsilon}}{2K} \, \Big| \,  {\cal T}_t^r  \right\} \\
& \leq & \sum_{k=1}^K \p \left\{ |W(t_k, \alpha_k) | \geq
\frac{t^\delta \, t_k^{1/\alpha_k}}{2K} \, \Big| \,  {\cal T}_t^r  \right\} \\
& = & \sum_{k=1}^K \p \left\{ |W(1, \alpha_k) | \geq
\frac{t^\delta}{2K} \right\},
\end{eqnarray*}
where we used that
$t^{v + \varepsilon} \geq t^\delta \,  t_1^{1/ \alpha_1}$ for some  small enough $\delta > 0$,
and that, by the definition of $A$ and $\varepsilon^\prime$, the inequalities
$t^{v + \varepsilon} \geq t^\delta t_j^{1/\alpha} \geq t^\delta t_j^{1/\alpha_j}$
hold, while in the last step the self-similarity of the stable process was used.
The last upper bound above goes to 0 as $t \to \infty$, and
$$
\sup_{t > 0} \E_i \sum_{j=1}^K N_t^j(\R^d) < \infty
$$
due to criticality of the branching. This finishes the proof of the lemma under the assumption
that $\gamma<1$. The proof for the case $\gamma = 1$  is a straightforward adaptation of the
previous one.
\end{proof*}

The value $\alpha_1 \gamma$ can be considered as the \emph{effective mobility} of
the type-1 particles. At an intuitive level if $\alpha_1 \gamma > \alpha$, then
second particle type is more mobile, even considering the long-living effect of
the first one, so that  in this case the ``dominant'' mobility is associated to
the second particle type. The next two theorems deal with the cases when the first
type is the dominant and when the second one, respectively. 

\begin{theorem}  \label{infmeanB-1dom}
Assume that \eqref{branchingassumption} holds and that $\gamma \eta> d/\alpha +1$.
If $ \alpha \geq \alpha_1 \gamma$,
i.e.  the mobility of the first particle type is dominant,
then the process suffers local extinction for $d < \alpha_1 \gamma / \beta$.
\end{theorem}

\begin{proof*}
Writing $u=t/2$ in Lemma \ref{L3VW}, and proceeding in the same way as we did to obtain
\eqref{q-ineq} in the proof of Theorem \ref{infmeanA}, we get
$$
q^{(i)}(t/2; 1 - c \, t^{-d/\alpha_1} ) \leq
c \, t^{-d/\alpha_1} t^{- \gamma \varepsilon/ (1 + \varepsilon)} +
c \, t^{- \gamma/(1 + \varepsilon) \beta}.
$$
Since in this case $v = 1 / \alpha_1$, from Lemma \ref{L13FVwith-v} we get extinction provided that
$$
\frac{d}{\alpha_1} < \frac{\gamma}{(1 + \varepsilon) \beta},
$$
which holds for $\varepsilon$ small enough if $d < \alpha_1 \gamma / \beta$.
\end{proof*}

\begin{theorem} \label{infmeanB-2dom}
Assume that $\gamma \eta > d/ \alpha + 1$.
If $\alpha_1 \gamma > \alpha$, i.e. the mobility of the second particle type is the dominant one, then
the process suffers local extinction for $ d < d_+$, where
\begin{equation} \label{d_+}
d_+ = \frac{ \gamma }
{\frac{( \beta + 1) \gamma}{\alpha} - \frac{1}{\alpha_1}}.
\end{equation}
\end{theorem}

\begin{proof*}
From the comparison lemma (Lemma \ref{comparison}) we have
$$
Q^{(i)}( t; \unit s ) \leq c\, n^{- \frac{1}{\beta}} + (1 - s)
\p \left\{ S_n \geq t \right\}.
$$
We have to  choose $t= n^{\frac{1 + \varepsilon}{\gamma}}$ for
some $\varepsilon > 0$, and then minimize the estimations in $\varepsilon$.
In this case
$$
q^{(i)}(t; 1 - s) \leq c\, t^{- \frac{\gamma}{(1+\varepsilon) \beta}} +
s \, t^{- \frac{\varepsilon \gamma}{1 + \varepsilon}}.
$$
Putting  $u=t/2$ in Lemma \ref{L3VWwith-nu} renders
$$
q^{(i)}(t/2; 1 - c_2 \, t^{-d/\alpha_1}) \leq c\, t^{- \frac{\gamma}{(1+\varepsilon) \beta}}
+ c \, t^{- d/\alpha_1 - \frac{\varepsilon \gamma}{1 + \varepsilon}}.
$$
Therefore we have to maximize
$$
\min \left\{ \frac{\gamma}{(1+ \varepsilon ) \beta } ,
\frac{d}{\alpha_1} + \frac{\varepsilon \gamma}{1+ \varepsilon} \right\}
$$
with respect to $\varepsilon$.
Since  the  term ${\gamma}/({(1+ \varepsilon ) \beta) }$ is monotone decreasing,
and the term $d / \alpha_1 + {\varepsilon \gamma}/({1+ \varepsilon})$
is increasing in $\varepsilon$, easy computations show that the optimal choice is
$$
\varepsilon = \frac{\gamma ( 1 + \beta^{-1})}{ d / \alpha_1 + \gamma} - 1,
$$
and the estimation is
$$
q^{(i)}(t/2; 1 - c_2 \, t^{- d / \alpha_1 }) \leq
c\, t^{- \frac{d / \alpha_1 + \gamma}{1 + \beta}}.
$$
Combining this with Lemma \ref{L13FVwith-v}, and taking into account
that $v= \gamma/ \alpha$, we get extinction if
$$
d \frac{\gamma}{\alpha} < \frac{ d / \alpha_1 + \gamma}{1 + \beta}.
$$
Solving the inequality, gives that extinction holds for $d < d_+$,
with the anticipated dimension $d_+$.
\end{proof*}\medskip

{\noindent\bf Remark\ }
 Notice that if $ \gamma / \alpha - 1 / \alpha_1 \to 0$, that is, if the effective mobilities
of types 1 and 2 are approximately the same, then $d_+ \to  \alpha_1 \gamma / \beta$, which
is the critical dimension in Theorem \ref{infmeanB-1dom}. Moreover,
for fixed $\alpha, \alpha_1$ and $\gamma$, the critical dimension $d_+$ considered as a
function of $\beta$, is decreasing, which is consistent with the known results.

\bigskip

\noindent
\textbf{Acknowledgement.}
We are grateful to the referee for the comments and remarks that greatly improved our paper.

\end{document}